\documentclass[11pt]{amsart}

\usepackage{amsmath}
\usepackage[font=small]{caption}
\usepackage{amsthm}
\usepackage{amsfonts}
\usepackage{graphicx}
\usepackage[margin=.5cm, format=hang]{subfig}
\usepackage[alphabetic]{amsrefs}
\usepackage{amssymb}
\usepackage{tikz-cd}
\usepackage{pinlabel}
\usepackage{morefloats}
\usepackage{enumitem}
\usepackage{fullpage}
\usepackage{array}
\usepackage{multirow}
\usepackage{sidecap}
\usepackage{hhline}
\usepackage{multirow}
\usepackage{blkarray}
\usepackage{xcolor}
\usepackage{comment}

\theoremstyle{plain}
\newtheorem{theorem}{Theorem}[section]
\newtheorem{proposition}[theorem]{Proposition}
\newtheorem{lemma}[theorem]{Lemma}
\newtheorem{corollary}[theorem]{Corollary}

\theoremstyle{definition}

\newtheorem{remark}[theorem]{Remark}

\usepackage{mathtools}

\newcommand{\R}{\mathbb{R}}

\newcommand{\ad}{\operatorname{ad}}
\author[Aliakbar Daemi]{Aliakbar Daemi}
\thanks{AD was partially supported by NSF Grants DMS-2609113, DMS-2208181, NSF FRG Grant DMS-1952762 and a Simons fellowship.}
\address{Department of Mathematics and Statistics, Washington University in St. Louis}
\email{adaemi@wustl.edu}

\author[Tye Lidman]{Tye Lidman}
\thanks{TL was partially supported by NSF Grants DMS-2105469 and DMS-2506277 and a Simons Travel Award.}
\address{Department of Mathematics, North Carolina State University}
\email{tlid@math.ncsu.edu}

\author[Mike Miller Eismeier]{Mike Miller Eismeier}
\thanks{MME was partially supported by NSTC grant 115-2115-M-001-011-MY3.}
\address{\vspace{-0.3cm}Institute of Mathematics, Academia Sinica}
\address{Department of Mathematics, University of Vermont}
\email{mme@as.edu.tw}
\email{Mike.Miller-Eismeier@uvm.edu}

\author[Nikolai Saveliev]{Nikolai Saveliev}
\thanks{NS was partially supported by the NSF FRG Grant DMS-1952762.}
\address{Department of Mathematics, University of Miami, Coral Gables, FL 33124}
\email{saveliev@math.miami.edu}

\title{Mazur manifolds and instanton Floer homology}

\begin{document}
\maketitle

\begin{abstract}
Instanton Floer homology of an integer homology 3-sphere is an 8-periodic homology theory that becomes 4-periodic with suitable coefficients. We prove that if an integer homology sphere bounds a Mazur manifold, then its instanton Floer homology is 2-periodic with any coefficients. Our proof relies on a variant of Floer's surgery exact triangle. We also obtain an obstruction to bounding a Mazur manifold in terms of certain values of the Chern--Simons invariants.
\end{abstract}

\section{Introduction}
In its original incarnation \cite{Floer1}, instanton Floer homology associates to an integer homology 3--sphere $Y$ an abelian group $I_*(Y)$. The group $I_*(Y)$, which is $\mathbb Z/8$ graded, is defined as the Morse homology group of the Chern--Simons functional of $Y$. As in other Floer homology theories, $I_*(Y)$ is generally well-behaved under cobordisms, and hence its structure has implications for the types of 4--manifolds that $Y$ can bound. For instance, half of the Euler characteristic of $I_*(Y)$ is equal modulo $2$ to the Rokhlin invariant $\mu(Y)$, which determines the signature modulo $16$ of any smooth, spin manifold with boundary $Y$; see \cite{AkbulutMcCarthy,Taubes-Casson}. Another example is Fr{\o}yshov's invariant $h(Y)$ defined using $I_*(Y)$ and certain auxiliary homomorphisms \cite{Froy1}. If $Y$ bounds a homology ball, the invariant $h(Y)$ must vanish; on the other hand, if a homology sphere $Y$ bounds a definite four-manifold with non-diagonalizable intersection form, $h(Y)$ must be nonzero; see \cite{Froy1, Froy2}. 

We add to this body of knowledge another result, which concerns the structure of $I_*(Y)$ in the case that $Y$ is the boundary of a Mazur manifold, a smooth $4$-manifold obtained by attaching a $2$-handle to $S^1 \times D^3$ along a curve which generates $H_1(S^1 \times S^2)$; see \cite{Mazur-AK}.  Such a manifold is necessarily contractible.  

\begin{theorem}\label{thm:mazur}
If $Y$ is the boundary of a Mazur manifold, then $I_k(Y;R) \cong I_{k+2}(Y;R)$ for all $k$ and all rings $R$.  
\end{theorem}

Let $R$ be a ring in which $2$ is invertible. Fr{\o}yshov proved that, for $k \ne 0,1 \pmod 8$, the $u$--map $ I_k(Y;R) \to I_{k-4}(Y;R)$ is an isomorphism \cite[Theorem 2]{Froy1}, so $I_*(Y;R)$ is 4--periodic. Froysh{\o}v's argument does not extend over $\mathbb Z$, and the computations in \cite{DMES} show that the instanton homology of $P \# P$, where $P = \Sigma(2,3,5)$ is the Poincar\'e homology sphere, is not $4$-periodic over either $\mathbb Z$ or $\mathbb Z/2$. 

The fact that Theorem~\ref{thm:mazur} holds over an arbitrary coefficient ring $R$ does carry interesting information. For example, $I_*(P \# -P; \mathbb Q) = 0$, according to \cite{Fukaya}. Therefore, the instanton Floer homology of $P \# -P$ over the rationals is $2$--periodic. However, 
\[
I_*(P \# -P; \mathbb Z) = (\mathbb Z/4)_{(0)} \oplus (\mathbb Z/2)^2_{(2)} \oplus (\mathbb Z/4)_{(4)} \oplus (\mathbb Z/2)^2_{(6)},
\]
according to the forthcoming work of Daemi--Miller Eismeier--Scaduto \cite{DMES}. Therefore, the instanton Floer homology of $P \# -P$ is not $2$--periodic over either $\mathbb Z$ or $\mathbb Z/2$. Theorem~\ref{thm:mazur} now implies that $P \# -P$ cannot bound a Mazur manifold.  This particular example, of course, can be deduced from other results. Akbulut's corollary of Taubes's periodic ends theorem implies that $P \# -P$ does not bound any positive-definite, simply-connected $4$-manifold \cite{Taubes}. Alternatively, the arguments of Gordon--Luecke \cite{GordonLuecke} imply that the boundary of any Mazur manifold is irreducible.

We now provide some historical context for Theorem~\ref{thm:mazur}.  By work of Galewski--Stern and Matumoto \cite{Tri-GS,Tri-Matumoto}, the high-dimensional triangulation conjecture is equivalent to the existence of an integer  homology sphere $Y$ such that $Y \# Y$ bounds a homology ball and $\mu(Y) = 1$. To disprove this, it suffices to construct a $\mathbb Z$--valued homology cobordism invariant which reduces modulo $2$ to $\mu$ and changes sign under orientation reversal.

Inspired by this, Saveliev \cite{Saveliev1} defined a $\mathbb{Z}$--valued invariant of integer homology spheres by 
\[
\nu(Y) = \frac{1}{2}\, \sum^7_{n=0} (-1)^{{(n+1)(n+2)}/{2}} \dim_{\mathbb{Q}} I_n(Y;\mathbb Q).  
\]
Modulo 2, the invariant $\nu$ is equal to half the Euler characteristic of $I_*(Y)$, and hence to the Rokhlin invariant. It also satisfies $\nu(-Y) = -\nu(Y)$. If $\nu$ were a homology cobordism invariant, this would disprove the triangulation conjecture. Saveliev \cite{Saveliev1, Saveliev2, Saveliev3} showed that $\nu$ vanishes for several classes of homology spheres known to bound acyclic manifolds. 

Much later, Manolescu completed the Galewski--Stern--Matumoto program by constructing, using $\text{Pin}(2)$-equivariant Seiberg--Witten theory, a $\mathbb Z$--valued homology cobordism invariant $\beta$ satisfying the conditions above; see \cite{Manolescu}. For Seifert fibered homology spheres, $\beta$ and $\nu$ coincide, although they do not agree in general \cite{Stoffregen}. Nevertheless, Theorem~\ref{thm:mazur} implies a new step toward Saveliev's original program:

\begin{corollary}\label{cor:nu}
If $Y$ is a homology sphere which bounds a Mazur manifold, then $\nu(Y) = 0$.  
\end{corollary} 

We now review the proof of Theorem~\ref{thm:mazur} and its mild enhancement. Suppose $Y$ is the boundary of a Mazur manifold, and let $Y'$ denote the homology sphere obtained by raising the framing of the 2--handle by 1. The instanton Floer homology groups of $Y$ and $Y'$ fit into an exact triangle \cite{Floer2} together with the manifold $S^1 \times S^2$ equipped with a nontrivial $SO(3)$--bundle. The latter manifold has trivial Floer homology, hence the map $I_*(Y) \to I_*(Y')$, induced by the cobordism map associated with a 2--handle attachment with trivial $SO(3)$--bundle, is an isomorphism of degree 0.

By modifying Scaduto's proof of the Floer exact triangle \cite{Scaduto}, we will show that $I_*(Y)$ and $I_*(Y')$ also fit into an exact triangle in which $S^1 \times S^2$ is equipped with the \textit{trivial} $SO(3)$--bundle. The resulting isomorphism $I_*(Y) \to I_*(Y')$ is induced by the same $2$--handle cobordism, now equipped with a nontrivial $SO(3)$--bundle, and is of degree $2$. 

Furthermore, we can easily determine the effect of these cobordism maps on Chern--Simons invariants. To this end, we again view instanton Floer homology heuristically as the Morse homology of the Chern--Simons functional. Just as a Morse function induces a filtration on its Morse complex, the Chern--Simons functional should induce a filtration on the instanton Floer chain complex. Making this precise requires both some subtlety with perturbations and a lift of the $\mathbb{Z}/8$--graded theory to a $\mathbb{Z}$--graded theory. Following, for example, the notation of \cite{DLME}, this yields an $\mathbb{R}$--filtered, $\mathbb{Z}$--graded abelian group $I_\bullet(Y)$. More precisely, $I_\bullet(Y)$ denotes the collection of groups $\{F_r I_*(Y)\}_{r \in \mathbb{R} \cup \infty}$, where $F_r I_*(Y)$ is, loosely speaking, the homology of the instanton Floer chain complex generated by flat connections with Chern--Simons invariant at most $r$. 

\begin{theorem}\label{thm:cs-drop}
Suppose that $Y$ is a homology sphere which bounds a Mazur manifold. Then there is a degree $-2$ automorphism of $I_\bullet(Y)$ which is strictly decreasing in the Chern--Simons filtration. There is also a degree $+2$ automorphism of $I_\bullet(Y)$ which increases the Chern--Simons filtration by strictly less than $1/2$. 
\end{theorem}  

The rest of the paper carries out the strategy outlined above to prove Theorems~\ref{thm:mazur} and \ref{thm:cs-drop}.  

\section*{Acknowledgements} We thank the Simons Center for Geometry and Physics for hosting the workshop ``Gauge theory and Floer homology in low dimensional topology" in 2025, where this project began. No generative AI was used in this project.

\section{The exact triangle}
Let $Y$ be a closed, oriented $3$-manifold. Recall that an oriented $1$-manifold $w \subset Y$ determines a $U(2)$-bundle $E_w \to Y$ with $c_1(E_w) = \text{PD}[w]$, unique up to essentially unique isomorphism; similarly, an oriented surface in an oriented $4$-manifold determines such a bundle \cite[Section 4.1]{KhUnknot}. If $(Y,w)$ is such a pair, its instanton Floer homology group $I_*(Y,w)$ is defined as follows \cite{Floer2}. Fix a connection $A_0$ on the determinant line bundle of $E_w$ and take the Morse homology of the Chern--Simons functional on the space of irreducible $U(2)$-connections of determinant equal to $A_0$, considered modulo determinant-$1$ gauge transformations. This is well-defined in favorable circumstances, including when $Y$ is an integer homology sphere or when $w$ has odd intersection number with some oriented surface in $Y$. Henceforth, we work with this Floer homology of $U(2)$--bundles. Up to isomorphism, $I_*(Y,w)$ only depends on the adjoint $SO(3)$-bundle and hence on the homology class of $w$ in $H_1(Y;\mathbb Z/2)$. 

For context, and to expand on the discussion in the introduction, recall that for a knot $K$ in a homology sphere $Y$, Floer \cite{Floer2,bd} constructed an exact triangle
\begin{equation}\label{eqn:Floer-triangle}
\ldots \to I_*(Y) \to I_*(Y_{-1}(K)) \to I_*(Y_0(K), w) \to \ldots 
\end{equation}
where $I_*(Y_0(K),w)$ is the instanton Floer homology of $Y_0(K)$ equipped with the $U(2)$--bundle dual to a generator $w$ of $H_1(Y_0(K))$. The maps in the exact triangle are induced by the natural 2--handle cobordisms, equipped with the corresponding $U(2)$--bundles: the cobordism $W: Y \to Y_{-1}(K)$ carries the trivial bundle, while the other two cobordisms carry non-trivial bundles.  In the special case when $Y_0(K) = S^1 \times S^2$, we have $I_* (Y_0(K),w) = 0$ and hence $W_*: I_*(Y) \to I_*(Y_{-1}(K))$ is an isomorphism of degree 0. 

The goal of this section is to prove the following result, which, when combined with the above discussion, establishes Theorem~\ref{thm:mazur}.  

\begin{proposition}\label{prop:deg-2-iso}
Suppose that $K$ is a knot in an integer homology sphere $Y$ and $Y_0(K) = S^1 \times S^2$. Equip the 2-handle cobordism $W: Y \to Y_{-1}(K)$ with the $U(2)$-bundle dual to the generator $c$ of $H_2(W, \partial W)$. For any coefficient ring, the induced map $(W,c)_*: I_*(Y) \to I_*(Y_{-1}(K))$ is an isomorphism of degree 2.
\end{proposition}

Proposition~\ref{prop:deg-2-iso} will follow from a new exact triangle, which we construct using standard gauge theoretic methods with a little additional care. We begin by recalling the general approach to exact triangles in instanton Floer homology following \cite[Section 5]{Scaduto}. Our only point of departure from that approach is the presence of reducible flat connections on $S^1 \times S^2$, which we address in Section~\ref{sec:reducibles}. Although not expressed in exactly this way, the argument constructs a homotopy inverse to the map $(W,c)_*$ in Proposition~\ref{prop:deg-2-iso}. The homotopy inverse is defined by counting ASD connections over a 1-parameter family of metrics on a certain cobordism from $Y_{-1}(K)$ to $Y$.  

\begin{remark}
The assumption that $Y_0(K) = S^1 \times S^2$ simplifies the analysis, but is not necessary for the existence of this exact triangle. To give a well-defined Floer homology of $Y_0(K)$, and establish an exact triangle as in Proposition~\ref{prop:deg-2-iso}, it is sufficient that the Alexander polynomial of $Y_0(K)$ have no roots on the unit circle.  However, we currently do not have any applications of this theory, and so do not develop this generalization here.  
\end{remark}

\subsection{The manifolds}\label{subsec:mfds}
Let $Y$ be an arbitrary closed oriented $3$-manifold equipped with a framed knot $K \subset Y$. Denote $Y_1 = Y$, $Y_2 = Y_{-1}(K)$, and $Y_0 = Y_0(K)$, and extend this notation to all $Y_i$ in a $3$--periodic fashion, $Y_i = Y_{i+3}$. Given a closed oriented $1$--manifold $w \subset Y \setminus K$, we equip the manifolds $Y_i$ with $U(2)$--bundles dual to  the respective geometric representatives $w_i$ given by
\[
w_1 = w_3 = w\quad\text{and}\quad w_2 = w \cup \widetilde K, 
\]
where $\widetilde K$ is the dual knot to the surgery. 

Next, we define a sequence of cobordisms $X_{ij}$ equipped with $U(2)$--bundles dual to geometric representatives $c_{ij} \subset X_{ij}$ with $c_{ij} \cap Y_i = w_i$ and $c_{ij} \cap Y_j = w_j$. The pairs $(X_{ij},c_{ij})$ will be defined for $1 \le j-i \le 3$ and extended to other indices in a $3$--periodic fashion, $X_{ij} = X_{i+3,j+3}$. The pair $(W,c)$ of Proposition~\ref{prop:deg-2-iso} will correspond to $(X_{12}, c_{12})$ with the empty $w$. 

Define $X_{i,i+1}$ to be the natural $2$--handle cobordism $I \times Y_i \cup D^2 \times D^2$ associated with the surgery on the respective knot. The cobordism $X_{12}$ is equipped with the bundle data $c_{12} = I \times w \cup cc_K$, where $cc_K$ is the co-core of the attached handle. The cobordism $X_{23}$ is equipped with the bundle data $c_{23} = I \times w \cup c_{\widetilde K}$, where $c_{\widetilde K}$ is the core of the handle attached along the dual knot $\widetilde K$. The cobordism $X_{34}$ is equipped with the bundle data $c_{34} = I \times w$. The remaining cobordisms $(X_{ij}, c_{ij})$ are obtained by stacking these basic cobordisms; for instance, 
\[
X_{14} = X_{12} \cup_{Y_2} X_{23} \cup_{Y_3} X_{34}.
\]

The construction of the exact triangle will rely on two observations about cobordisms $X_{ij}$. The first one is that the cobordism $X_{i-1,i+1}$ contains a canonical $(-1)$-sphere, given as the union of the co-core of $X_{i-1,i}$ with the core of $X_{i,i+1}$. Write $S_i$ for the boundary of a tubular neighborhood of this $(-1)$-sphere; this separates $X_{i-1,i+1}$ as a union 
\[
X_{i-1,i+1} = X'_{i-1,i+1} \cup_{S_i} U,\quad\text{where}\;\; U = \overline{\mathbb{CP}}^2 \setminus B^4.
\] 
Crucially, the bundle data $c_{i-1,i+1} \cap U$ gives rise to the non-trivial $SO(3)$-bundle over $U$, as in \cite[Section 3.4]{Scaduto}. This can be understood on a case-by-case basis. In the case of $X_{13}$, this is because $c_{13}$ is the canonical $(-1)$-sphere in $U$. In the case of $X_{02}$ and $X_{24}$, this is because $c_{02}$ and $c_{24}$, respectively, have intersection number $1$ with the canonical $(-1)$-spheres in these manifolds.

\begin{figure}[!ht]
\labellist
\pinlabel $Y_1$ at 80 520
\pinlabel $Y_2$ at 305 520
\pinlabel ${Y_0}$ at 550 520
\pinlabel $Y_1$ at 750 520
\pinlabel $S_0$ at 490  350
\pinlabel $S_2$ at 350 350
\pinlabel $T$ at 630 375
\endlabellist
\centering
		\includegraphics[width=0.7\textwidth]{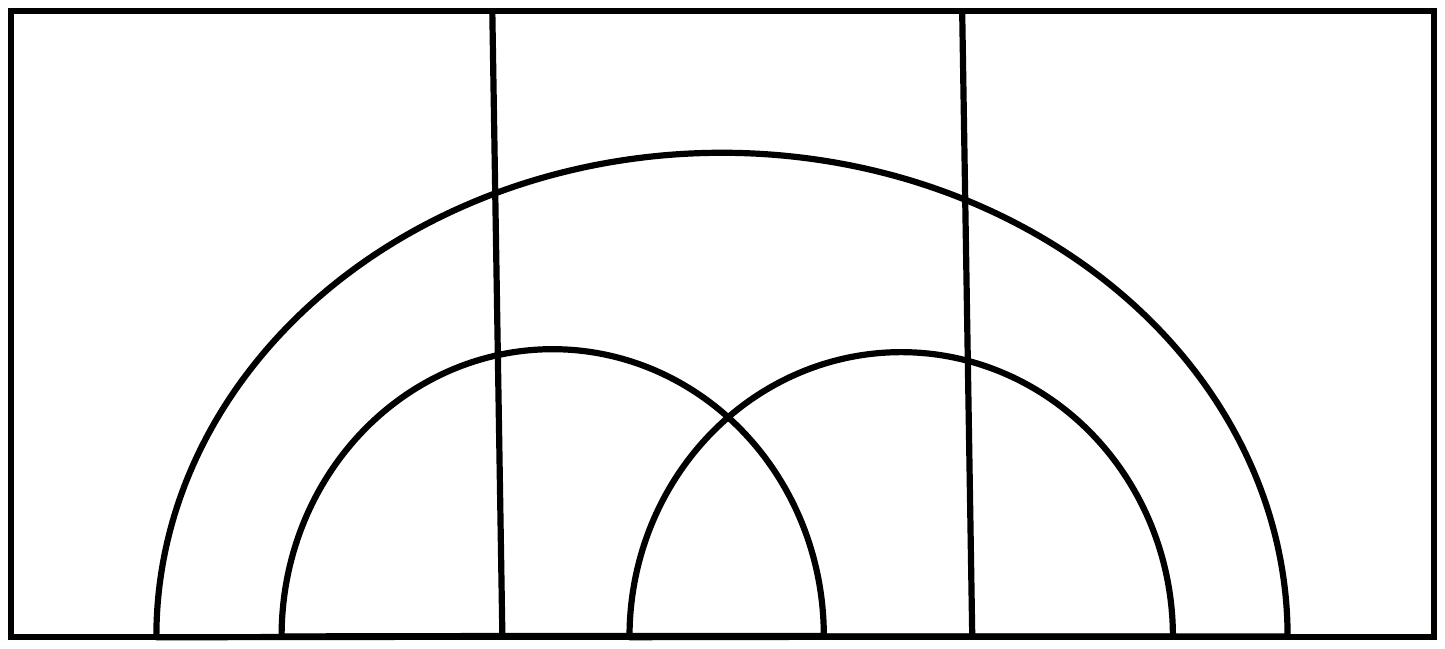}
        $ $\vspace{-.85in}
		\caption{The cobordism $X_{14}: Y_1 \to Y_1$ with the five canonical hypersurfaces $Y_0$, $Y_2$, $S_0$, $S_2$, and $T$.}
		\label{F3}
\end{figure}

The other observation uses the fact that the inclusions of $X_{i,i+2}$ and $X_{i+1,i+3}$ give rise to the $(-1)$-spheres $S_{i+1}$ and $S_{i+2}$ in $X_{i,i+3}$, which intersect transversely at a single point. Let $V$ be a regular neighborhood of the union of these spheres. Then we have a decomposition 
\[
X_{i,i+3} = X'_{i,i+3} \cup_T V,\quad\text{where}\;\; V = (D^2 \times S^2)\, \#\, \overline{\mathbb{CP}}^2\;\;\text{and}\;\; T = \partial V = S^1 \times S^2.
\]
The manifold $X'_{i,i+3}$ here is diffeomorphic to the complement of a tubular neighborhood of $\{1/2\} \times K \subset I \times Y_i$. Note that the $SO(3)$ bundle over $T$ is trivial. The cobordism $X_{i,i+3}$ is equipped with five canonical hypersurfaces: $Y_{i+1}$, $Y_{i+2}$, $S_{i+1}$, $S_{i+2}$, and $T$. In the case $i = 1$, these are depicted in Figure \ref{F3}. 

Fix Riemannian metrics on each $Y_i$, and equip each $X_{i,i+1}$ with a choice of metric, cylindrical in a chosen collar neighborhood of the boundary corresponding to the fixed metrics on $Y_i$, $Y_{i+1}$. This allows us to consider the corresponding noncompact, complete Riemannian manifold obtained by attaching cylindrical ends to the boundary of $X_{i,i+1}$; by an abuse of notation, we still call this manifold $X_{i,i+1}$. 

Owing to the two hypersurfaces $Y_i$ and $S_i$, we may choose a family of metrics $G^{i-1}_{i+1}$ on $X_{i-1,i+1}$ parameterized by an interval. As we approach one endpoint of the interval, the metrics develop an increasingly long neck stretched along $Y_i$; as we approach the other endpoint, they develop an increasingly long neck stretched along $S_i$. The boundary points of the interval correspond to `broken' or `fully-stretched' metrics, the first equal to the disjoint union of $X_{i-1,i}$ and $X_{i,i+1}$, and the latter equal to the disjoint union of $X'_{i-1,i+1}$ and $U$, all equipped with cylindrical-end metrics. 

In a similar fashion $X_{i,i+3}$ admits a family of metrics $G^i_{i+3}$ parameterized by a pentagon; metrics near the boundary of the pentagon have a long neck along at least one of the five canonical hypersurfaces. The five boundary faces correspond to metrics which are fully stretched along one of the five hypersurfaces, while the five vertices correspond to metrics which are fully stretched along a pair of disjoint hypersurfaces. 

\subsection{Triangle detection}\label{sec:triangle}
One effective way to prove the existence of an exact triangle is the {\em triangle detection lemma}, see e.g. \cite[Lemma 4.2]{OzSz:HF-br}, \cite[Lemma 7.1]{KhUnknot} or \cite[Lemma 5.1]{Scaduto}. 

\begin{lemma}[Triangle detection lemma]\label{lem:triangle-detection}
Suppose $(C_i, d_i)$ are chain complexes, with $(C_{i+3},d_{i+3}) = (C_i, d_i)$. Suppose further that we have maps \[f_i: C_i \to C_{i+1}, \quad h_i: C_i \to C_{i+2}, \quad k_i: C_i \to C_{i+3}\] satisfying the following conditions: 
\begin{align*}
d_{i+1} f_i - f_i d_i &= 0, \\
d_{i+2} h_i + h_i d_i + f_{i+1} f_i &= 0, \\
d_{i+3} k_i - k_i d_i - f_{i+2} h_i - h_{i+1} f_i &= q_i,
\end{align*}
where $q_i$ induce isomorphisms on homology. Then the $f_i$ define an exact triangle: \[\cdots \to H(C_1, d_1) \xrightarrow{f_1} H(C_2, d_2) \xrightarrow{f_2} H(C_3, d_3) \xrightarrow{f_3} H(C_1, d_1) \to \cdots\]
\end{lemma}

We will apply this lemma in the case when $(C_i,d_i)$ is the Floer chain complex $CI(Y_i,w_i)$, whose generators are irreducible flat connections on the bundle $E_{w_i} \to Y_i$ modulo gauge equivalence, with differential $d_i$ counting irreducible ASD connections of index $1$ on the pullback bundle of $E_{w_i}$ to $\R \times Y_i$, considered modulo translation and gauge equivalence. For this to make sense, $E_{w_i}$ should support a finite number of irreducible flat connections, all cut out transversely. Similarly, the moduli spaces of irreducible ASD connections on $\R \times Y_i$ of index at most 2 should be cut out transversely, and the moduli space of index 1 should be finite after factoring out the translation symmetry. This can be arranged with a suitable choice of perturbations. Further perturbations of the ASD equation on various cobordisms $(X_{ij},c_{ij})$, which agree with the given perturbations on the ends, may be needed to make the relevant moduli spaces cut out transversely. Our treatment of perturbations is identical to that in \cite{Scaduto}, so we will omit further mention of them in the rest of this section.

The maps $f_i: C_i \to C_{i+1}$ are defined by counting irreducible ASD connections of index $0$ on $(X_{i,i+1}, c_{i,i+1})$, with irreducible flat limits on $Y_i$ and $Y_{i+1}$. The maps $h_i: C_i \to C_{i+2}$ are defined by counting irreducible ASD connections of index $-1$ on $(X_{i,i+2}, c_{i,i+2})$ with irreducible flat limits, parameterized by the family $G^i_{i+2}$. Finally, the maps $k_i: C_i \to C_{i+3}$ are defined by counting irreducible ASD connections of index $-2$ on $(X_{i,i+3}, c_{i,i+3})$ with irreducible flat limits, parameterized by the family $G^i_{i+3}$. Note that all these parameterized moduli spaces of ASD connections have dimension zero.

A version of the argument of \cite[Section 5]{Scaduto} now applies to establish the hypotheses of the triangle detection lemma. Similar arguments appear in \cite{KhUnknot,KMOS:len-surgery,PFH}.

\begin{itemize}
	\item The fact that $f_i$ is a chain map follows by counting the ends of the 1-dimensional moduli space of irreducible connections of index $1$. 
	\item The chain-homotopy relation for $h_i$ follows from counting the ends of the 1-dimensional moduli space of irreducible connections on $(X_{i,i+2}, c_{i,i+2})$ of index 0 parameterized by $G^i_{i+2}$. In addition to energy sliding off the ends, this moduli space also has ends corresponding to limiting towards the broken metrics. The metric fully stretched along $Y_{i+1}$ contributes $f_{i+1} f_i$. There are no further contributions because metrics with long necks along $S_i$ carry no irreducible ASD connections of index $\le 0$. This follows as in \cite[Section 5.2]{Scaduto} from the additivity of the index and the determination of the smallest-index ASD connections on both sides of $S_i$; here it is crucial that $U$ is equipped with the non-trivial $SO(3)$-bundle. 
	\item The relation on $k_i$ follows from counting the ends of the 1-dimensional moduli space of irreducible connections on $(X_{i,i+3}, c_{i,i+3})$ of index $-1$ parametrized by $G^i_{i+3}$. Energy sliding off the ends gives rise to $d_{i+3} k_i - k_i d_i$, while ASD connections supported by metrics fully stretched along $Y_{i+1}$ or $Y_{i+2}$ give rise to $-f_{i+2} h_i$ and $-h_{i+1} f_i$. There are no ASD connections of index $-1$ with a long neck along $S_{i+1}$ or $S_{i+2}$ by the same index-additivity argument as in the previous bullet point. Finally, the metric fully-stretched along $T$ gives the term $m$, which is chain homotopic to an isomorphism by arguing that $(V, G^i_{i+3})$ can be replaced by $S^1 \times D^3$ with trivial bundle data and a single standard metric; replacing $V$ with $S^1 \times D^3$ returns $Y_i \times I$.  
\end{itemize}

Compactness, which is needed for these arguments to make sense, would follow directly from the Uhlenbeck compactness theorem, were it not for the possible appearance of reducible connections. A sequence of ASD connections limits to a broken ASD connection supported by a broken metric, and the limiting connection may in principle have reducible components. The argument can tolerate reducible connections on $U$ and $V$, but not on $X_{ij}$ or $X'_{ij}$, and it cannot tolerate reducible flat limits on any $Y_i$.

This is our point of departure from \cite{Scaduto}. There, the author rules out reducible ASD connections by requiring each $(Y_i, w_i)$ to be \textit{admissible}, meaning that $[w_i] \in H_1(Y_i;\mathbb Z/2)$ pairs oddly with some orientable surface. It then follows that $(Y_i, w_i)$ supports no reducible flat connections, and hence that neither $X_{ij}$ nor $X'_{ij}$ can support a reducible ASD connection.  To establish the exact triangle asserted in Proposition \ref{prop:deg-2-iso}, we must give a similar argument for the manifolds in question, based on a careful analysis of their reducible ASD connections.  More specifically, for the cobordisms appearing in Floer's original exact triangle \eqref{eqn:Floer-triangle}, the only reducible connections occuring in the moduli spaces above are those broken along $T$, as described in the last bullet point. We must prove the analogous statement in our setting. 

\subsection{The case of interest}
We will now restrict ourselves to $Y_0 = S^1 \times S^2$ and empty $w$. The manifolds $Y_1$ and $Y_2$ are integer homology spheres hence the bundles over them can be trivialized using a trivialization which is unique up to homotopy, so this data is inconsequential. The trivial bundle over $Y_0 = S^1 \times S^2$ carries no irreducible flat connections hence $C_0 = 0$ and the maps $f_0$, $f_2$, $h_0$, $h_1$ and $k_0$ are automatically zero.  All that is needed is the following technical result.

\begin{proposition}\label{prop:index}	
Let the manifolds $Y_i$ and $X_{jk}$ be as above and assume that all the required transversality conditions hold. For $jk \in \{12,24,14\}$, let $A$ be a broken ASD connection on $X_{jk}$ such that 
\begin{equation}\label{E:hypothesis}
i(A) + \dim G^j_k\; \le\; 1,
\end{equation}
and suppose that $A$ has irreducible flat limits on $Y_j$ and $Y_k$. Then either $A$ consists entirely of irreducible ASD connections with irreducible flat limits, or else $jk = 14$ and $A$ is broken along $T$, with an irreducible ASD connection on the $X'_{14}$ component.
\end{proposition}

We will now complete the proof of the main results of the paper, assuming Proposition~\ref{prop:index}, which will be proved in the next and final section.

\begin{proof}[Proof of Proposition \ref{prop:deg-2-iso}]
Because $f_0$ is zero and Proposition \ref{prop:index} implies that the relevant relations in the triangle detection lemma hold, it follows that $h_2$ is a chain map and the composition $h_2 f_1$ of the maps $f_1: I_*(Y_1) \to I_*(Y_2)$ and $h_2: I_*(Y_2) \to I_*(Y_1)$ is an isomorphism. Because Floer's exact triangle implies that $I_*(Y_1)$ and $I_*(Y_2)$ are isomorphic, and surjections of isomorphic finitely-generated modules are isomorphisms, $f_1$ and $h_2$ are in fact isomorphisms.
\end{proof}

\begin{proof}[Proof of Theorem \ref{thm:cs-drop}]
We begin by observing that the cobordisms $X_{i,i+1}$ and hence all cobordisms $X_{ij}$ are simply-connected. This can be seen as follows. Since $Y_1$ and $Y_2$ are obtained by integral surgery along a knot of winding number one in $S^1 \times S^2$, their fundamental groups are normally generated by the dual knot.  Attaching a 2-handle along that knot kills the fundamental group of both $X_{01}$ and $X_{23}$. The composite of $X_{23}$ and $X_{01}$ is the blowup of the orientation-reversal of $X_{12}$, hence $X_{12}$ is also simply connected.
\end{proof}

We now turn to the degree and the behavior with respect to the Chern--Simons functional. A degree $D$ endomorphism of $I_*(Y)$ is of {\it level} $L$ if it maps $F_r I_d(Y)$ into $F_{r+L} I_{d+D}(Y)$ for any $r$ and $d$. 
The map $f_1$ is induced by the simply-connected cobordism $X_{12}$ with $b^+(X_{12}) = 0$ and $c_{12}^2 = -1$, so the formula of \cite[Proposition 3.4]{DLME} gives that $f_1$ has degree $2$ and level less than $1/4$. The map $h_2$ is induced by a $1$-parameter family $G^2_4$ of metrics on a simply-connected cobordism $X_{24}$ with $b^+(X_{24}) = 1$ and equipped with a trivializable bundle, so a similar argument shows that $h_2$ has degree $-2$ and negative level. In Floer's original exact triangle, the map analogous to $f_1$ has degree $0$ and negative level, while the map analogous to $h_2$ has degree $0$ and level less than $1/4$. Composing the appropriate maps gives an automorphism of $I_*(Y_1)$ of degree $2$ and level less than $1/2$ and another automorphism of degree $-2$ and negative level.

\section{Reducibles in the exact triangle}\label{sec:reducibles}
In this section, we will prove Proposition~\ref{prop:index}. Much of the background material can be found in \cite{DonBook}. To set conventions, a connection on an $SU(2)$-bundle is called \textit{central} if its holonomy lies in the center of $SU(2)$, \textit{abelian} if its holonomy group is abelian but non-central, and \textit{irreducible} if its holonomy group is nonabelian. 

Our goal is to show that the broken ASD connections of Proposition~\ref{prop:index} consist of irreducible ASD connections with irreducible limits, with the one noted exception. We will use families of perturbations of the ASD equation which are trivial on the ends corresponding to $S^3$ or $S^1 \times S^2$, and equal to a small fixed perturbation over the ends corresponding to $Y_1$ or $Y_2$, chosen so that all unbroken irreducible ASD connections over a family of metrics $G$, with respect to the chosen family of perturbations, are cut out transversely. Provided that the perturbations are sufficiently small, the central instantons will be cut out transversely.  We will omit further mention of perturbations. The transversality condition ensures that, for any irreducible ASD connection $A$ supported by the family of metrics $G$, we have
\begin{equation}\label{E:trans}
0\; \le\; i(A) + \dim G,
\end{equation}
where $i(A)$ is the index of $A$, that is, the index of the respective ASD operator with the Atiyah--Patodi--Singer boundary conditions \cite{APS1}.

\subsection{Some index theory}
Let $W$ be a 4-manifold with a principal $U(2)$-bundle $E$ and ends $M_1,\ldots,M_k$.  Let $A$ be a connection on $E$ which limits to the flat connections $\alpha_1,\ldots,\alpha_k$ on the ends. When $b_1(W) = 0$, the index of $A$ is shown in \cite[Equation 82]{MMR} to be
\begin{equation}\label{eqn:index-formula}
i(A) = 8 \kappa(A) - \frac{3}{2}(\chi(W) + \sigma(W)) - \frac 1 2\, \sum\, h(\alpha_{\ell}) + \frac 1 2\,\sum\, \rho(\alpha_{\ell}).
\end{equation}

The terms in this formula are as follows. For a flat $SU(2)$ connection $\alpha$ on a 3-manifold $Y$, $h(\alpha) = h^0(\alpha) + h^1(\alpha)$, where $h^j(\alpha) = \dim H^j(Y; \R^3_{\alpha})$ and $\R^3_{\alpha}$ is the local system with fiber $\R^3$ and holonomy $\ad\alpha$. Note that, on a homology sphere, the unique central flat connection has $h(\theta) = 3$, while non-degenerate flat irreducibles have $h(\alpha) = 0$. On $S^1 \times S^2$, the two central connections have $h(\theta) = 6$, while each abelian connection has $h(\alpha) = 2$. 

The term $\rho(\alpha)$ is the Atiyah--Patodi--Singer $\rho$-invariant of a connection $\ad\alpha$ on the adjoint $SO(3)$-bundle over $Y$ with the boundary orientation; see \cite{APS2}. If $\alpha$ is a central connection, or any flat connection on $S^1 \times S^2$, then $\rho(\alpha) = 0$. 

The term $\kappa(A)$ is the topological energy of $A$, defined as 
\[
\kappa(A) = -\frac 1 4\, p_1 (\ad A) = \frac{1}{32\pi^2} \displaystyle\int \operatorname{Tr}( F (\ad A) \wedge F(\ad A)),
\] 
where $F(\ad A)$ is the curvature of the adjoint connection $\ad A$. Recall that, for an ASD connection $A$, we have $\kappa (A) \geq 0$, with equality if and only if $A$ is projectively flat. In addition, if the $M_{\ell}$ are homology spheres or $S^1 \times S^2$ and the $\alpha_{\ell}$ are all reducible, then $8\kappa(A) \in 2\mathbb Z$. This is because the manifold $W$ can be closed up by attaching any 4-manifold to $M_{\ell}$ if $\alpha_\ell$ is the trivial connection, and attaching $S^1 \times D^3$ to $S^1 \times S^2$ with a flat reducible connection extending the connection on $S^1 \times S^2$. Since the integral is zero on all added pieces, $8 \kappa(A)$ is $-2$ times the first Pontryagin number on a closed manifold, hence twice an integer. If $A$ is an $SU(2)$--connection, $p_1 (\ad A) = c_1 (A)^2 - 4 c_2 (A) = - 4 c_2 (A)$ so $8 \kappa(A) \in 8\mathbb Z$ by a similar argument. 

Next, suppose that  $W_1$, $W_2$ are $4$-manifolds with cylindrical ends, with $W_1$ containing an end $[0,\infty) \times M$ and $W_2$ containing an end $(-\infty, 0] \times M$. If $A_1$ and $A_2$ are connections on $W_1$ and $W_2$ which are asymptotic to a common flat connection $\alpha$ on $M$, then we can glue $W_1$ and $W_2$ together along this end and obtain a new connection $A_1 * A_2$. The index formula \eqref{eqn:index-formula} gives
\[
i(A_1 * A_2) = i(A_1) + i(A_2) + h(\alpha).   
\]

It is convenient to consider an alternative quantity that is additive under gluing: for an asymptotically flat connection $A$ on $W$, define
\[
\tau(A) = i(A) + \frac 1 2\,\sum h(\alpha_{\ell}). 
\]  
Note that, when the limits of $A$ are non-degenerate irreducibles, $\tau(A) = i(A)$, and when the limits of $A$ are reducibles on a homology sphere or $S^1 \times S^2$, 
\[
\tau(A) = 8\kappa (A) - \frac{3}{2}\, (\chi(W) + \sigma(W)).
\]

\subsection{Estimates on $\tau (A)$}
The transversality condition ensures that non-constant ASD connections on a cylinder $\R \times Y$ are irreducible and have index at least $1$. It follows that if $A$ is a non-constant ASD connection on a cylinder equipped with the trivial bundle and limits $\alpha, \beta$ then the following lower bounds hold:
\begin{itemize}
\item If both $\alpha, \beta$ are irreducible, then $\tau(A) \geq 1$; 
\item If exactly one of $\alpha, \beta$ is reducible, then $\tau(A) \geq 5/2$;
\item If both of $\alpha, \beta$ are reducible, then $\tau(A) \geq 8$. 
\end{itemize}
We would also like to have bounds on $\tau(A)$ for ASD connections $A$ on the $4$-manifolds discussed in Section \ref{subsec:mfds}. These are summarized in Table \ref{tab:r}. For simplicity, we suppress the bundle data, although it is important to keep in mind that the bundle is non-trivial on several of these manifolds.  

\medskip

\begin{table}[!ht]
\begin{tabular}{ |p{1.5cm}|p{1.2cm}|p{2.3cm}|p{2cm}|p{1.1cm}|p{2cm}|  }
\hline
 Manifold & $\chi + \sigma$ & $8\kappa \geq$ (reduc) & $\tau \geq$ (reduc) & $\dim G $ & $\tau \geq$ (irred)\\
 \hline
 $X_{12}$	&0	 	&2 	&2		&0	&0	 \\
 $X_{23}$ 	&1 		&0	&$-3/2$	&0	&1        \\
 $X_{34}$ 	&1 		&0 	&$-3/2$	&0	&1	 \\
 $X_{24}$	&2 		&2	&$-1$ 	&1	&$-1$     \\
 $X'_{24}$	&1		&0	&$-3/2$	&0	&3/2	\\
 $U$    	&1 		&2 	&1/2	&0	&3/2     \\
 $X_{13}$ 	&1  	&2	&1/2	&1	&0     \\
 $X'_{13}$ 	&0 		&0 	&0 		&0	&5/2      \\
 $X_{14}$ 	&2  	&2	&$-1$  	&2 	&$-2$     \\
 $X'_{14}$ 	&0  	&0 	&0	&0 	&1 \\
 $V$	&2 		&2	&$-1$   	&1	&0    \\
% $V'$ 	&1 		&0	&$-3/2$	&0	&5/2     \\
 \hline
\end{tabular}
\caption{Lower bounds on $\tau(A)$ for ASD connections $A$ on relevant 4-manifolds}
\label{tab:r}
\end{table}

The fourth column of Table \ref{tab:r} records the least possible value of $\tau(A)$ of an ASD connection $A$ with reducible flat limits; $A$ itself need not be reducible. The final column records the least possible value of $\tau(A)$ of an irreducible ASD connection $A$ supported by the relevant family of metrics $G$. By our transversality hypothesis \eqref{E:hypothesis}, any irreducible ASD connection $A$ (possibly broken, but with no reducible components) satisfies
\[
\tau(A) \geq -\dim G + \frac 1 2\,\sum h(\alpha_{\ell}) \geq - \dim G.
\] 
The first inequality is stronger when one of the ends corresponds to $S^1 \times S^2$ or $S^3$. 

The computation of $\chi + \sigma$ in Table \ref{tab:r} is purely topological and is omitted. The lower bound on $8 \kappa (A)$ for reducible connections follows from the fact that the ends correspond to homology spheres or $S^1 \times S^2$. The third column records whether the bundle admits a projectively flat connection. The fourth column is then simply $8\kappa(A) - 3(\chi + \sigma)/2$. 

\subsection{Analysis of broken ASD connections} 
To prove Proposition~\ref{prop:index}, we need to analyze the possibly broken ASD connections on the various $X_{ij}$, for relevant metrics, and determine when they may fail to be irreducible or to have irreducible limits. The proof of the following two lemmas is identical to that of \cite[Lemma 5.5]{DLME}.
  
\begin{lemma}\label{lem:fully-irreducible}
If $A$ is a broken ASD connection as in Proposition~\ref{prop:index}, and all components of $A$ are irreducible, then all flat limits of $A$ are irreducible.
\end{lemma}
\begin{proof}
Write $A = A_1 * A_2$, where the $A_j$ are possibly broken irreducible ASD connections supported by families of metrics $G_j$ with a common flat limit $\alpha$. Then either one of the $A_j$ is a nonconstant ASD connection on a cylinder supported by a single metric and the other is supported by $G$, or else $G_1 \times G_2$ is a boundary face of $G$, so that $\dim G = \dim G_1 + \dim G_2 + 1$. In the first case, $i(A_1) + i(A_2) \ge 1 - \dim G$, since the nonconstant ASD connection on the cylinder has index at least 1 and the other connection has index bounded below by $-\dim(G)$ by hypothesis \eqref{E:trans}. In the second case, hypothesis \eqref{E:trans} gives $i(A_1) + i(A_2) \ge -\dim G_1 - \dim G_2 = 1 - \dim G$. In either case, 
\[
i(A) = i(A_1) + i(A_2) + h(\alpha) \ge 1-\dim G + h(\alpha).
\]
Now, hypothesis \eqref{E:hypothesis} requires $i(A) \le 1 - \dim G$, hence $h(\alpha) = 0$ and $\alpha$ must be irreducible.
\end{proof}

\begin{lemma}\label{lem:no-red-breaking}
If $A$ is a broken ASD connection on $X_{jk}$ as in Proposition~\ref{prop:index}, then $A$ does not break along $Y_j$ or $Y_k$ with a reducible limit. 
\end{lemma}
\begin{proof}
Write $A = A_j * B * A_k$, where $A_j$ and $A_k$ are possibly empty broken trajectories on $\mathbb R \times Y_j$ and $\mathbb R \times Y_k$, respectively, and $B$ is a broken ASD connection on $X_{jk}$ without any components supported on the cylinders $\R \times Y_j$ or $\R \times Y_k$. If both outer limits of $B$ are reducible, we have $\tau(B) \ge -3/2$ by Table \ref{tab:r}. Because $A$ has irreducible outer limits, we also have $\tau(A_j) \ge 5/2$ and $\tau(A_k) \ge 5/2$ using estimates on the cylinder. Therefore, $\tau(A) \ge 7/2$, which contradicts hypothesis \eqref{E:hypothesis}. If $B$ has at least one irreducible limit, $B$ must be irreducible, so we have $\tau(B) \ge i(B) \ge -\dim G^j_k$ by \eqref{E:trans}. A non-constant ASD connection on a product $\mathbb R \times Y$ must have $\tau(A) \ge 1$, hence it follows from hypothesis \eqref{E:hypothesis} that at most one of $A_j$, $A_k$ can be nonempty and, if nonempty, it must consist of an irreducible trajectory between irreducible flat connections.
\end{proof}

\begin{proof}[Proof of Proposition~\ref{prop:index}]
By Lemma \ref{lem:fully-irreducible}, it suffices to show that $A$ has no reducible components, except when $jk= 14$, with $A$ broken along $T$. By Lemma \ref{lem:no-red-breaking}, any reducible component has neither $Y_j$ nor $Y_k$ its boundary component. Furthermore, reducible ASD connections on cylinders are constant, and $A$ includes no constant trajectories. We call a non-cylindrical codimension-$0$ submanifold with neither $Y_j$ nor $Y_k$ in its boundary an \textit{interior piece}; we aim to show $A$ has no reducible component on an interior piece. We consider the three cobordisms in turn:

\medskip

\begin{itemize}
\item{\em Cobordism $X_{12}$.} No broken configuration has an interior piece.
 \end{itemize}
 
 \smallskip
 
 \begin{itemize}
 \item{\em Cobordism $X_{24}$.} The only potential configuration $A$ with an interior piece would be one broken along $S_0$. We will show that no such configuration occurs. Indeed, if $A$ were broken along $S_0$, then its $X'_{24}$ component would satisfy $\tau \ge 3/2$, while its $U$ component would satisfy $\tau \ge 1/2$. Hence $\tau(A) \ge 2$, contradicting hypothesis \eqref{E:hypothesis}, which gives $\tau(A) \le 0$.
\end{itemize}

\smallskip

\begin{itemize}
\item {\em Cobordism $X_{14}$.} A broken metric on $X_{14}$ with an interior piece is broken along $S_0, S_2$, or $T$, or $T$ together with one of the $S_i$. Configurations broken only along $S_0$ or $S_2$ are ruled out by the same analysis as for $X_{24}$. If $A$ is reducible on $X_{23}$, then $\tau(A) \ge 0 - 3/2 + 1 > -2$. Thus $A$ must be broken along $T$, with its restriction to $X'_{14}$ irreducible. \qedhere
\end{itemize}
\end{proof}
\smallskip

\begin{remark}
It is straightforward to verify that our proof of Proposition \ref{prop:index} holds for all $j,k$ with $\dim G^j_k \le 2$, giving a slightly more direct proof of Proposition \ref{prop:deg-2-iso}. Similarly, it is possible to verify that Proposition \ref{prop:index} holds for $S^1 \times S^2$ replaced by any manifold with first homology $\mathbb Z$ and Alexander polynomial with no roots on the unit circle.
\end{remark}

\bibliography{paper}

\end{document}